\documentclass[10pt,twoside]{article}
\usepackage[latin1]{inputenc}
\usepackage{amsmath}
\usepackage{graphicx}
\usepackage{yhmath}
\usepackage[table]{xcolor}
\usepackage{ mathrsfs} 
\usepackage{amssymb}
\usepackage{url}

\input epsf
\def\limiten{\renewcommand{\arraystretch}{0.5}
\begin{array}[t]{c}\stackrel{}{\longrightarrow} \\
{\scriptstyle n\rightarrow
\infty}\end{array}\renewcommand{\arraystretch}{1}}

\def\limitepsn{\renewcommand{\arraystretch}{0.5}
\begin{array}[t]{c}\stackrel{a.s.}{\longrightarrow} \\
{\scriptstyle n \rightarrow
\infty}\end{array}\renewcommand{\arraystretch}{1}}

\DeclareMathOperator\argmax{argmax}

\renewcommand{\Box}{\hfill\rule{0.25cm}{0.25cm}} 

\newtheorem{Prop}{Proposition}[section]
\newtheorem{lem}{Lemma}[section]
\newtheorem{Theo}{Theorem}[section]
\newtheorem{cor}{Corollary}[section]
\newtheorem{rem}{Remark}[section]

\newcommand{\E}{\ensuremath{\mathbb{E}}}
\newcommand{\R}{\ensuremath{\mathbb{R}}}
\newcommand{\Z}{\ensuremath{\mathbb{Z}}}

\newcommand{\N}{\ensuremath{\mathbb{N}}}

\definecolor{grisclair}{gray}{0.9}

\newcommand*\interior[1]{\overset{\mathsf{o}}{#1}}

\renewcommand{\arraystretch}{.8}

\begin{document}
\title{\bf Strong consistent model selection for general causal time series}
 \maketitle \vspace{-1.0cm}
 \begin{center}
   William Kengne \footnote{Developed within the ANR BREAKRISK (ANR-17-CE26-0001-01) }. 
 \end{center}

  \begin{center}
 CY Cergy Paris Université, CNRS, THEMA, F-95000 Cergy, France.\\
{\it  E-mail:  william.kengne@u-cergy.fr }
 \end{center}

 \pagestyle{myheadings}
 \markboth{Strong consistent model selection}{W. Kengne}

\textbf{Abstract} :
 We consider the strongly consistent question for model selection in a large class of causal time series models, including  AR($\infty$), ARCH($\infty$), TARCH($\infty$), ARMA-GARCH and many classical others processes.
 We propose a  penalized criterion based on the quasi likelihood of the model.
 We provide sufficient conditions that ensure the strong consistency of the proposed procedure. Also, the estimator of the parameter of the selected model obeys the law of iterated logarithm. 
 It appears that, unlike the result of the weak consistency obtained by Bardet {\it et al.} \cite{Bardet2020}, a dependence between the regularization parameter and the model structure is not needed.  \\


 \section{Introduction}
 We consider a general class of autoregressive time series in a semiparametric framework. 
  Let $M , f : \R^{\N} \rightarrow \R$ be two measurable functions and $(\xi_t)_{t\in \Z}$ a sequence of centered independent
 and identically distributed (iid)  random variables satisfying $\textrm{var}(\xi_0)= 1$.
 %
 Consider the class of affine causal models,

\medskip

\noindent \textbf{Class} $\mathcal{AC}(M,f):$ A process $X=(X_t)_{t\in \Z}$ belongs to $\mathcal{AC}(M,f)$ if it satisfies:
\begin{equation}
X_t=M\big((X_{t-i})_{i\in \mathbb{N}^*}\big)\, \xi_t+f\big((X_{t-i})_{i\in \mathbb{N}^*}\big) \;\; \mbox{for any}~ t\in \Z.
\label{model}
\end{equation}
The existence of a stationary and ergodic solution of the class (\ref{model}) has been studied by Bardet and Wintenberger (2009) as a particular case of models considered in Doukhan and Wintenberger (2008). Bardet and Wintenberger (2009) and Bardet {\it et al.} (2017) carried out the inference question in the semiparametric setting in the class $\mathcal{AC}(M,f)$, whereas Bardet {\it et al.} (2012), Kengne (2012), Bardet and Kengne (2014) focussed on the change-point problem in this class.   
Numerous classical time series models belongs to the class (\ref{model}): for instance AR$(\infty)$, ARCH$(\infty)$, TARCH$(\infty)$, ARMA-GARCH or APARCH processes. 
 
 \medskip

 Consider a trajectory $(X_1,\ldots,X_n)$ of a process $X=(X_t)_{t\in \Z}$ that belongs to $\mathcal{AC}(M^*,f^*)$  where $M^*$ and $f^*$ are unknown.
We consider a finite collection ${\cal M}$ of affine causal models, where the true model $m^* \in \cal M$ corresponds to $M^*$ and $f^*$.
 Our main aim is to select a model $\widehat m$ (among the collection $\cal M$) which is "close" to $m^*$ for large $n$.

 \medskip

 We focus on the semiparametric framework and assume that the distribution of $\xi_0$ is unknown and that the functions $f$ and $M$ are known up to a parameter $\theta \in \Theta$, where $\Theta$ is a compact subset of $\R^{d}$ ($d \in \N$). That is, the model $m^*$ corresponds to the true parameter $\theta^* \in \Theta$ and the process $X$ belongs to  $\mathcal{AC}(M_{\theta^*},f_{\theta^*})$.
 In the sequel :
 \begin{itemize}
\item each model $m \in \cal M$ is considered as a subset of $\{1,\ldots,d\}$ and denote by $|m|$ the dimension of the model, typically, $|m|=\#(m)$;
\item for $m \in \cal M$, the parameter space of $m$ is $\Theta(m)=\big \{ (\theta_i)_{1\leq i \leq d}\in \Theta,~\theta_i=0~\mbox{if $i\notin m$} \big \}$; $\theta(m)$ is the parameter vector associated to $m$;
\item the collection ${\cal M}$ is considered as a subset of the power set of $\{1,\ldots,d\}$,  i.e. ${\cal M}\subset {\cal P} \big ( \{1,\ldots,d\}\big )$.
\end{itemize}
Therefore, for any model $m\in {\cal M}$, $m \in \mathcal{AC}(M_{\theta},f_{\theta})$ when $\theta \in \Theta(m)$.
Also, we could consider hierarchical as well as exhaustive families of models. 
 
 \medskip

\noindent  For instance, assume that $(X_1,\ldots,X_n)$ is generated from a GARCH($p^*,q^*$) process. The collection ${\cal M}$ could be a family of ARMA($p,q$)-GARCH($p',q'$) with $(p,q,p',q') \in \{0,1,\cdots,p_{max} \} \times \{0,1,\cdots,q_{max} \} \times \{0,1,\cdots,p_{max}' \} \times \{0,1,\cdots,p_{max}' \}$ where $p_{max}, q_{max}, p_{max}', p_{max}'$  are the fixed upper bounds of the orders, assumed to satisfy $ p^* \leq p_{max}$ and $ q^* \leq q_{max}$. 
Therefore, consider $\Theta$ as a compact subset of $\R^{p_{max} + q_{max}} \times (0,\infty) \times [0,\infty)^{p_{max}' + q_{max}'}$. Thus, a model $m$ is a subset of $\{1,\cdots, p_{max} + q_{max} + p_{max}' + q_{max}' +1\}$ and its parameter space is $\Theta(m)=\big \{ (\theta_i)_{1\leq i \leq d}\in \Theta,~\theta_i=0~\mbox{if $i\notin m$} \big \}$. 
 
 \medskip

%

 The model selection problem for time series has already been considered by several authors; we refer to the book of McQuarrie and Tsai (1998), the monograph of Rao and Wu (2001), the recent review paper of Ding {\it et al.} (2018), the recent works of Hsu {\it et al.} (2019) and the references therein for an overview on this topic.
 Hannan (1980) and Hannan and Deistler (2012) provided general conditions for strong consistency of the order estimator of an ARMA and ARMAX model. 
 Resende and Dorea (2016) proposed the efficient determination criterion   (introduced by Zhao {\it et al.} (2001) for the strongly consistent estimation of the order of multiple Markov chains) for model selection in a general class of multivariate time series. 
 They established the strong consistency of the procedure under some conditions which may seem a bit strong for some applications, for instance, the existence of the third order derivative of the contrast function (likelihood), the existence of moments of order 16 for the BEKK-GARCH model.
 Recently, Bardet {\it et al.} (2020) addressed the model selection question in the class of model $\mathcal{AC}(M_{\theta},f_{\theta})$. They proposed a procedure based on the quasi likelihood of the model and provided sufficient conditions that ensure the weak consistency of the selected model. 
  
 \medskip
 
 In this new contribution, we focus on the model selection in the class of model $\mathcal{AC}(M_{\theta},f_{\theta})$ with a penalized contrast which is based on the Gaussian quasi likelihood of the model.
 \begin{itemize}
  \item[(i)] Under the assumptions that $\E|\xi_0|^r < \infty$ with $r>4$, the functions $\theta \mapsto f_\theta, ~ M_\theta$ are twice times continuously differentiable on $\Theta$ and satisfy some Lipschitz-type properties, we establish the strong consistency of the proposed procedure.
  \item[(ii)] We show that the quasi maximum likelihood estimator (QMLE) of the selected model obeys the law of iterated logarithm.
  \end{itemize} 

 \medskip
 
 The rest of the paper is structured as follows. In Section 2, we set some notations, assumptions and define the model selection criterion.  The main results are provided in Section 3 whereas Section 4 is devoted to a concluding remarks. Section 5 focuses on the proofs of the main results. 
 

\section{Assumptions and the model selection criterion}
\subsection{Assumptions on the class of models $\mathcal{AC}(M_\theta,f_\theta)$}
 In the sequel, we will use the norms :
 \begin{enumerate}
 \item  $\|\cdot\|$ applied to a vector denotes the Euclidean norm of
    the vector;
    \item for any compact set $ \mathcal{K} \subseteq\R^d$ and for any
    $g:\Theta \longrightarrow\R^{p}$, $ \|g\|_{ \mathcal{K}}=\sup_{\theta\in \Theta}(\|g(\theta)\|)$.
 \end{enumerate}
 Throughout the sequel, we assume that the functions $ \theta \mapsto M_\theta$ and $\theta \mapsto f_\theta$
 are twice times continuously differentiable on $\Theta$.
Also, we will use $H_{\theta} = M_{\theta}^2$ and for any function $g_{\theta}$ which is $i$ times  differentiable on $\Theta$, we set $\partial_{\theta^i} g_{\theta} = \partial^i g_{\theta} / \partial\theta^i$.
Let us consider the following assumptions for any compact set $\mathcal{K} \subseteq \Theta$, $i=0,1,2$ and $\Psi_{\theta}=\partial_{\theta^i} f_{\theta}$ or $\partial_{\theta^i}M_{\theta}$ :

\medskip

\noindent \textbf{Assumption A}$(\Psi_{\theta},\mathcal{K})$: for any $x \in \R^\N$, the function $\theta \mapsto \Psi_{\theta}(x)$ is continuous on $\Theta$ with $\|\Psi_{\theta}(0)\|_\Theta < \infty$ and there exists a sequence of non-negative real numbers $\big(\alpha_k(\Psi_{\theta},\mathcal{K})\big)_{k\ge 1}$ satisfying  $\sum_{k=1}^{\infty}\alpha_k(\Psi_{\theta},\mathcal{K})< \infty$  such that:
	\[
	\|\Psi_{\theta}(x)-\Psi_{\theta}(y)\|_\mathcal{K} \le \sum_{k=1}^{\infty}\alpha_k(\Psi_{\theta},\mathcal{K})|x_k-y_k| \; for \; all\; x,y \in \mathbb{R}^\N.
	\]
 \noindent In the sequel we refer to the particular case of the non linear ARCH($\infty$) (NLARCH$(\infty)$, see Bardet and Wintenberger(2009)) processes define when $f_\theta =0$. In this case,  we consider the following assumption  for $i=0,1,2$ :

\medskip

\noindent \textbf{Assumption A}$(\partial_{\theta^i} H_\theta,\mathcal{K})$: for any $x \in \R^\infty$, the function $\theta \mapsto \partial_{\theta^i} H_\theta(x)$ is continuous on $\Theta$ with $\|\partial_{\theta^i} H_\theta(0)\|_\Theta < \infty$ and there exists a sequence of non-negative real numbers $\big(\alpha_k(\partial_{\theta^i} H_\theta,\mathcal{K})\big)_{k\ge 1}$ satisfying $\sum_{k=1}^{\infty}\alpha_k(\partial_{\theta^i} H_\theta,\mathcal{K})< \infty$ such that :
	\[
	\|\partial_{\theta^i} H_\theta(x)-\partial_{\theta^i} H_\theta(y)\|_\mathcal{K} \le \sum_{k=1}^{\infty}\alpha_k(\partial_{\theta^i} H_\theta,\mathcal{K})|x_k^2-y_k^2| \; for \; all\; x,y \in \mathbb{R}^\N.
	\]

\medskip

Then define the set:
  \begin{multline}\label{def_Theta_r}
  \Theta(r)=\Big \{\theta  \in \R^d,~ A(f_{\theta},\{\theta\})\; \textnormal{and}\; A(M_{\theta},\{\theta\})\; \textnormal{hold with} 
\sum_{k=1}^{\infty} \alpha_k(f_{\theta},\{\theta\}) +\|\xi_0\|_r \, \sum_{k=1}^{\infty} \alpha_k(M_{\theta},\{\theta\}) < 1 \Big \} \\
\cup 
\Big \{\theta  \in \R^d,~ f_\theta = 0\; \textnormal{and}\; A(H_{\theta},\{\theta\})\; \textnormal{holds with } \|\xi_0\|_r^2 \, \sum_{k=1}^{\infty} \alpha_k(H_{\theta},\{\theta\}) < 1 \Big \} 
 \end{multline}
 The above Lipschitz-type conditions are classical when studying the existence of a stationary and ergodic solution of such model, see Doukhan and Wintenberger (2008). In the case of the class $\mathcal{AC}(M_\theta,f_\theta)$, if $\theta\in\Theta(r)$, then there exists a unique causal, stationary and ergodic  solution  $X=(X_t)_{t\in \Z}\in \mathcal{AC}(M_\theta,f_\theta)$ with a finite moment of order $r$, see Bardet and Wintenberger (2009).

\medskip

The following assumptions are useful in the study of the asymptotic behavior of the QLME .\\
 \noindent \textbf{Assumption D}$(\Theta)$: {\it $\exists \underline{h}>0$ such that $\underset{\theta \in \Theta}{\inf}(H_{\theta}(x)) \ge \underline{h}$ for all $x\in \mathbb{R}^{\N}$.} 
 
\medskip

\noindent \textbf{Assumption Id}$(\Theta)$:  For all  $(\theta,\theta')\in \Theta^2$,
\[ \Big( f_{\theta}(X_0,X_{-1},\cdots)=f_{\theta'}(X_0,X_{-1},\cdots)~\text{and}~M_{\theta}(X_0,X_{-1},\cdots)=M_{\theta'}(X_0,X_{-1},\cdots) \ \text{a.s.}\Big) \Rightarrow \ \theta = \theta'.\]

\medskip

\noindent  { \textbf{Assumption Var$(\Theta)$}: For all $\theta  \in \Theta $, one
 of the families $ \big( \dfrac{\partial f_{\theta}}{\partial \theta^{i}}(X_0,X_{-1},\cdots) \big)_{1\leq i \leq d} \quad
  \mbox{or} \quad \big( \dfrac{\partial h_{\theta}}{\partial \theta^{i}}(X_0,X_{-1},\cdots) \big)_{1\leq i \leq d}  \quad $
 is a.s. linearly independent.
 
 \medskip

    \noindent In the following assumption, we make the convention that if \textbf{A}$(\partial_{\theta^i} M_\theta,\Theta)$ holds then $\alpha_k(\partial_{\theta^i} H_\theta,\Theta) = 0$ and if
     \textbf{A}$(\partial_{\theta^i} H_\theta,\Theta)$ holds, then $\alpha_k(\partial_{\theta^i} M_\theta,\Theta) = 0.$ Set for $\ell = 0, 1, 2$,

\medskip
 \noindent \textbf{Assumption $\textbf{K}_\ell(\Theta)$}: there exists $r > 4$ such that $\theta^* \in \Theta(r) \cap \Theta$ and  for $i=0,\cdots,\ell$,  $A(\partial_{\theta^i} f_\theta,\Theta), A(\partial_{\theta^i} M_\theta,\Theta)$ (or $A(\partial_{\theta^i} H_\theta,\Theta)$) hold  with
\[
\sum_{k\ge 1} \frac{1}{ \sqrt{k \log \log k} } \sum_{j \ge k} \sum_{i=0}^\ell  \alpha_j (\partial_{\theta^i} f_\theta,\Theta)+ \alpha_j (\partial_{\theta^i} M_\theta,\Theta)+\alpha_j (\partial_{\theta^i} H_\theta,\Theta)  < \infty.
\] 
These aforementioned assumptions hold for many classical models, including AR($\infty$), ARCH($\infty$), TARCH($\infty$) type processes , see for instance Bardet and Wintenberger (2009),  Bardet {\it et al.} (2012), Kengne (2012).
In the case of assumption $\textbf{K}(\Theta)$, let us consider for $\ell=0, 1, 2$ :
\begin{enumerate}
\item The geometric case: $ \sum_{i=0}^\ell  \alpha_j (\partial_{\theta^i} f_\theta,\Theta)+ \alpha_j (\partial_{\theta^i} M_\theta,\Theta)+\alpha_j (\partial_{\theta^i} H_\theta,\Theta) = \mathcal{O}(a^j)$ for some $a \in [0,1)$. In this case, assumption $K_\ell(\Theta)$ holds.
\item The Riemanian case:  $ \sum_{i=0}^\ell  \alpha_j (\partial_{\theta^i} f_\theta,\Theta)+ \alpha_j (\partial_{\theta^i} M_\theta,\Theta)+\alpha_j (\partial_{\theta^i} H_\theta,\Theta) = \mathcal{O}(j^\gamma)$ with $\gamma >0$. If $\gamma > 3/2$, then $K_\ell(\Theta)$ holds. 
\end{enumerate}

 \subsection{The model selection criterion}
 Consider a model $m \in {\cal M}$ and the class $\mathcal{AC}(M_{\theta},f_{\theta})$ for $\theta \in \Theta(m) \subset \Theta \subset \R^d$.
 Assume that a trajectory $(X_1,\ldots,X_n)$ is observed. 
The conditional Gaussian quasi (log)likelihood (up to a constant) $L_n$ is defined for all $\theta \in \Theta(m)$ by,
\begin{equation}
L_n(\theta):=-\frac{1}{2}\sum_{t=1}^n q_t (\theta) ~ , ~ \textnormal{with} \; q_t (\theta):=\frac{(X_t - f_{\theta}^t)^2}{H_{\theta}^t} + \log(H_{\theta}^t)
\label{eq:eq1}
\end{equation}
where $f_{\theta}^t:=f_{\theta}(X_{t-1},X_{t-2},\cdots)$, $M_{\theta}^t:=M_{\theta}(X_{t-1},X_{t-2},\cdots)$ and $H_{\theta}^t=\big (M_{\theta}^t\big )^2$.
Since $L_n(\theta)$ depends on $(X_t)_{t\le 0}$ which are not observed, it is common practice (see \cite{Bardet2009}, \cite{Bardet2012}, \cite{Kengne2012}) to consider the approximated  quasi (log)likelihood  given (up to a constant)  for all $\theta \in \Theta(m)$ by
\begin{equation}\label{eq:eq2}
\widehat{L}_n(\theta):=-\frac{1}{2}\, \sum_{t=1}^n \widehat{q}_t (\theta) ~ , ~ \textnormal{with} \; \widehat{q}_t (\theta):=\frac{(X_t - \widehat{f}_{\theta}^t)^2}{\widehat{H}_{\theta}^t} + \log(\widehat{H}_{\theta}^t) 
\end{equation}
 where $\widehat{f}_{\theta}^t=f_{\theta}(X_{t-1},X_{t-2},\cdots,X_1,0,\cdots)$, $\widehat{M}_{\theta}^t=M_{\theta}(X_{t-1},X_{t-2},\cdots,X_1,0,\cdots)$, $\widehat{H}_{\theta}^t=(\widehat{M}_{\theta}^t)^2$.
 Note that, the "best" parameter associated to the model $m$ if defined by,
 \[ \theta^*(m) = \underset{\theta \in \Theta(m)}{\text{argmin} ~\E[q_0(\theta)]}  .\]
According to \cite{Bardet2020}, $\theta^*(m)$ exists and it is unique under Id$(\Theta(m))$. When $m=m^*$, we have $\theta^*(m^*)=\theta^*$. 
For any $m \in \cal M$, the QMLE of $\theta^*(m)$ is given by 
\begin{equation}
\widehat{\theta}(m)= \underset{\theta \in \Theta(m)}{\argmax} \;\widehat{L}_n(\theta).
\label{eq:qmle}
\end{equation}
The selection of the "best" model $\widehat{m}$  among the collection $\mathcal{M}$ is performed by minimizing the penalized contrast
\begin{equation}
\widehat{C}(m)=-2\widehat{L}_n\big(\widehat{\theta}(m)\big)+ |m| \kappa_n,
\label{eq:cri}
\end{equation}
that is
\begin{equation}
\widehat{m}=  \underset{m \in \mathcal{M}}{\text{argmin}} ~\widehat{C}(m),
\label{est_m}
\end{equation}
where
\begin{itemize}
	\item $(\kappa_n)_n$ is the sequence of the regularization parameter  (possibly data-dependent) that will be used to calibrate the penalty term;
	\item $|m|$ is the dimension of the model $m$, typically, the cardinal of $m$ (considered as a subset of $\{1,\ldots,d\}$), which is also the number of the estimated components of $\theta$ (the others are fixed to zero).
\end{itemize}

\section{Asymptotic results}
 Recall that, when the model is correctly specified, Bardet and Wintenberger (2009) have established the consistency and the  asymptotic normality of $\widehat\theta(m^*)$. The following proposition shows that the estimator $\widehat\theta(m^*)$ obeys the law of iterated logarithm. 
\begin{Prop}\label{proposition}
Let $(X_1,\ldots, X_n)$ be a trajectory of a process $X$ belonging to $\mathcal{AC}(M_{\theta^*},f_{\theta^*})$ 
where $\theta^* \in \Theta(r) \cap \interior{\Theta } \subset \R^d$ with $r> 4$.  Assume that $D(\Theta)$, $Id(\Theta)$, $Var(\Theta)$, $K_2(\Theta)$  hold. Then,
\begin{equation} 
\widehat\theta(m^*) - \theta^* = \mathcal{O} \Big( \sqrt{\dfrac{\log \log n}{n}} \Big) ~ ~ a.s. ~ .
\label{prop_eq}
\end{equation}
\end{Prop}

\medskip

\noindent The following theorem provides sufficient conditions that ensure the strong consistency of the model selection procedure.
%
\begin{Theo}\label{theorem}
Let $(X_1,\ldots, X_n)$ be a  trajectory of a process $X$ belonging to $\mathcal{AC}(M_{\theta^*},f_{\theta^*})$). 
 Under the assumptions of Proposition \ref{proposition}, and if $\kappa_n/\log \log n \limiten \infty$ and $\kappa_n/n \limiten 0$, then
\begin{equation}
\widehat{m} \limitepsn m^*  .
\label{theo_eq}
\end{equation}
\end{Theo}

\medskip

\begin{rem}
\begin{enumerate}
 \item This result, besides it is stronger than those obtained by Bardet {\it et al.} (2020), do not impose any condition on the dependence between the regularization parameter $\kappa_n$ and the Lipschitz-type coefficients $\alpha_j (\partial_{\theta^i} f_\theta,\Theta)$, $\alpha_j (\partial_{\theta^i} M_\theta,\Theta)$, $\alpha_j (\partial_{\theta^i} H_\theta,\Theta)$ as have been set by these authors. For instance, in the Riemanian case with $ \sum_{i=0}^\ell  \alpha_j (\partial_{\theta^i} f_\theta,\Theta)+ \alpha_j (\partial_{\theta^i} M_\theta,\Theta)+\alpha_j (\partial_{\theta^i} H_\theta,\Theta) = \mathcal{O}(j^\gamma)$ for some $3/2<\gamma <2$, the BIC ($\kappa_n= \log n$) is strongly consistent from this Theorem while the result of Bardet {\it et al.} (2020) can not assure the weak consistency of the BIC. 
 \item Hannan and Deistler (2012) have considered the estimation of the order of an ARMAX (including ARMA), where the contrast $\widehat{C}$ is based on the  Gaussian likelihood of the model. 
 Under the condition $\kappa_n/n \limiten 0$, they have established that there exists a constant $c_1 >0$ such that, if $\liminf \kappa_n / (2\log \log n) > c_1$, then the estimator of the order is strongly consistent.
 From the proof of Theorem \ref{theorem}, one can see that such result holds for the general class of model considered here; that is, we can find a constant $c_2>0$ such that if $\liminf \kappa_n / \log \log n > c_2$, then $\widehat{m} \limitepsn m^*$.
 \end{enumerate} 

\end{rem}

\medskip

 \noindent The next corollary show that the estimator of the parameter of the selected model  $\widehat\theta(\widehat{m})$ obeys the law of iterated logarithm.
\begin{cor}\label{corollary}
Let $(X_1,\ldots, X_n)$ be a trajectory of a process $X$ belonging to $\mathcal{AC}(M_{\theta^*},f_{\theta^*})$. 
 Under the assumptions of Theorem \ref{theorem},
\begin{equation} 
\widehat\theta(\widehat{m}) - \theta^* = \mathcal{O} \Big( \sqrt{\dfrac{\log \log n}{n}} \Big) ~ ~ a.s. ~ .
\label{cor_eq}
\end{equation}
\end{cor}

\section{Concluding remarks} 
This paper focuses on the model selection in a large class of causal time series models in a semiparametric framework. 
 The strong consistency of an estimator based on a penalized quasi likelihood contrast is established, under some classical conditions on the regularization parameters $\kappa_n$. \\
 \indent For the estimation of the order of an ARMAX model, Hannan and Deistler (2012) have established that, there exists a constant $c_0$ such that, if $\limsup \kappa_n / (2\log \log n) < c_0$ then the strong consistency of the estimator of the order fails. A topic of a future works could be to investigate if such result is applied to the general class of model considered here or to derive an upper bound of $\kappa_n$ for which the strong consistency fails. \\
 Another extension of this works is to carry out the model selection problem in the class $\mathcal{AC}(M_{\theta},f_{\theta})$ with a procedure based on a non Gaussian (for instance Laplacian, see Bardet {\it et al.} (2017)) quasi likelihood, for the purpose of reducing the order moment imposed on the process.

\section{Proofs of the main results}\label{Proofs}
The following lemma will be useful in the sequel.
The proof is carried out by going along similar lines as in Lemma 2 of \cite{Bardet2020} by using  Corollary 1 of \cite{Kounias1969}; so, it is then omitted. 

\begin{lem} \label{lem1}
Let $X \in \mathcal{AC}(M_{\theta},f_{\theta})$  and $\Theta \subseteq \Theta(r)$  with {$r > 4$}. Assume that the assumptions $D(\Theta)$ and $K_1(\Theta)$  hold. Then,
\begin{equation}\label{eq_lem1}
\frac{1}{ \sqrt{n \log \log n} } \,\Big \|\dfrac{\partial\widehat{L}_n(\theta)}{\partial \theta}-\frac{\partial L_n(\theta)}{\partial \theta}\Big \|_{\Theta}\limitepsn 0.
\end{equation}
\end{lem}
%
%
  
\medskip

\paragraph{Proof of Proposition \ref{proposition}}~ 
According to \cite{Bardet2009}, it holds that  $ \widehat{\theta}(m^*) \limitepsn \theta^*$. Also, since $\theta^* \in \Theta(m^*) \cap \interior{\Theta}$, we get $ \dfrac{\partial\widehat{L}_n(\widehat{\theta}(m^*) )}{\partial \theta} =0 $ for $n$ large enough. 
Thus, for any $i=1,\cdots,|m^*|$, the Taylor expansion of  $\dfrac{\partial\widehat{L}_n}{\partial \theta_i}$ implies
\[ 0 =  \dfrac{\partial\widehat{L}_n(\widehat{\theta}(m^*) )}{\partial \theta_i} = \dfrac{\partial\widehat{L}_n(\theta^*)}{\partial \theta_i} +   \dfrac{\partial^2 \widehat{L}_n\big( \dot{\theta}_i(m^*) \big)}{\partial \theta \partial \theta_i} \big( \widehat{\theta}(m^*) - \theta^* \big), \]
where $\dot{\theta}_i(m^*)$ lies between $\widehat{\theta}(m^*)$ and $\theta^*$.
Therefore,
\begin{equation}\label{Proof_prop_theta_m_theta_star_partial_widehat_L}
\sqrt{\dfrac{n}{\log \log n}}\big( \widehat{\theta}(m^*) - \theta^* \big) = \dfrac{2}{ \sqrt{n \log \log n}} \widehat{F}_n^{-1}(m^*) \dfrac{\partial\widehat{L}_n(\theta^*)}{\partial \theta} ~ ~ \text{ where } ~ 
 \widehat{F}_n(m^*) = -2 \Big( \dfrac{\partial^2 \widehat{L}_n\big( \dot{\theta}_i(m^*) \big)}{\partial \theta \partial \theta_i} \Big)_{i \in m^*}. 
\end{equation} 
Note that, by dealing with the first (stationary) regime in the Corollary 6.1 of \cite{Bardet2012} and since $\dot{\theta}_i(m^*) \limitepsn \theta^*$ for $i=1,\cdots,|m^*|$, we get
\begin{equation}\label{Proof_prop_widehat_F_conv}
\widehat{F}_n(m^*) \limitepsn F(\theta^*,m^*) ~  \text{ where } ~ F(\theta^*,m^*)= \Big(\E \Big[ \dfrac{ \partial^2 q_0(\theta^*)}{ \partial\theta_i \partial\theta_j}  \Big] \Big)_{i,j \in m^*} . 
\end{equation}   
Since $F(\theta^*,m^*)$ is invertible (see \cite{Bardet2009}), then for $n$ large enough and with a sufficiently large probability, the matrix $\widehat{F}_n(m^*)$ is invertible.   
We have from Lemme \ref{lem1}, (\ref{Proof_prop_theta_m_theta_star_partial_widehat_L}) and (\ref{Proof_prop_widehat_F_conv}), 
\begin{equation}\label{Proof_prop_theta_m_theta_star_partial_L}
\sqrt{\dfrac{n}{\log \log n}}\big( \widehat{\theta}(m^*) - \theta^* \big) = \dfrac{2}{ \sqrt{n \log \log n}} \widehat{F}_n^{-1}(m^*) \dfrac{\partial L_n(\theta^*)}{\partial \theta} + o(1) ~ ~ a.s. ~ . 
\end{equation} 
 We have,
\[  \dfrac{\partial L_n(\theta^*)}{\partial \theta} = \sum_{t=1}^n \dfrac{\partial q_t(\theta^*)}{ \partial \theta} .\]
Denote for all $t \in \Z$, $\mathcal{F}_t = \sigma(X_t,X_{t-1},\cdots)$ the $\sigma$-field generated by the whole past at time $t$.  Then, $\Big( \frac{\partial q_t(\theta^*)}{ \partial \theta},\mathcal{F}_t \Big)$ is a stationary ergodic square integrable martingale difference process (see \cite{Bardet2009}). Therefore, from the law of iterative logarithm for martingales (see \cite{Stout1970, Stout1974}), we get,
\begin{equation}\label{Proof_prop_LIL}
\frac{1}{ \sqrt{n \log \log n}}  \dfrac{\partial L_n(\theta^*)}{\partial \theta} = \mathcal{O}(1) ~ ~ ~ a.s. ~.
\end{equation}
Thus, the proposition follows from (\ref{Proof_prop_widehat_F_conv}), (\ref{Proof_prop_theta_m_theta_star_partial_L}) and (\ref{Proof_prop_LIL}).
 \begin{flushright}
$\blacksquare$ 
\end{flushright}

\paragraph{Proof of Theorem \ref{theorem}}~ 

\medskip

1. Let $m\in {\cal M}$ such as $m  \supsetneq  m^*$. 
We have, 
%
\begin{equation}\label{C_m_star_C_m_log_log_n}
 \frac{1}{ \log \log n} \big( \widehat{C}(m^*) - \widehat{C}(m) \big) = \frac{2}{\log \log n} \big(\widehat{L}_n\big(\widehat{\theta}(m)\big)-\widehat{L}_n\big(\widehat{\theta}(m^*)\big) - \frac{\kappa_n}{\log \log n}  ( |m|-|m^*|). 
\end{equation}
Let us establish that
\begin{equation}\label{L_m_star_L_m_log_log_n}
 \frac{1}{\log \log n} \big(\widehat{L}_n\big(\widehat{\theta}(m)\big)-\widehat{L}_n\big(\widehat{\theta}(m^*)\big) = \mathcal{O}(1) ~ a.s. ~.
\end{equation}
Since $\theta^* \in \Theta(m) \cap \interior{\Theta}$ and $ \widehat{\theta}(m) \limitepsn \theta^*$, then $ \dfrac{\partial\widehat{L}_n(\widehat{\theta}(m) )}{\partial \theta} =0 $ for $n$ large enough. 
Therefore, from the Taylor expansion of $\widehat{L}_n$, we can find $\overline{\theta}(m)$ between $\widehat{\theta}(m)$ and $\theta^*$ such that
\begin{equation}\label{Taylor_widehat_L}
\widehat{L}_n\big(\widehat{\theta}(m)\big)-\widehat{L}_n\big(\theta^*) =  \frac{1}{2}\big( \widehat{\theta}(m) - \theta^* \big)' \dfrac{\partial^2 \widehat{L}_n\big( \overline{\theta}(m) \big)}{\partial \theta^2} \big( \widehat{\theta}(m) - \theta^* \big).
\end{equation}
Also, for any $i=1,\cdots,|m|$, we can find $\dot{\theta}_i(m)$ between $\widehat{\theta}(m)$ and $\theta^*$ such that, for $n$ large enough,
 \[ 0 =  \dfrac{\partial\widehat{L}_n(\widehat{\theta}(m) )}{\partial \theta_i} = \dfrac{\partial\widehat{L}_n(\theta^*)}{\partial \theta_i} +   \dfrac{\partial^2 \widehat{L}_n\big( \dot{\theta}_i(m) \big)}{\partial \theta \partial \theta_i} \big( \widehat{\theta}(m) - \theta^* \big) .\]
Hence,
\begin{equation}\label{theta_m_theta_star_partial_widehat_L}
\widehat{\theta}(m) - \theta^* = \dfrac{2}{n} \widehat{F}_n^{-1}(m) \dfrac{\partial\widehat{L}_n(\theta^*)}{\partial \theta} ~ ~ \text{ where } ~ 
 \widehat{F}_n(m) = -2 \Big( \dfrac{\partial^2 \widehat{L}_n\big( \dot{\theta}_i(m) \big)}{\partial \theta \partial \theta_i} \Big)_{i \in m}. 
\end{equation}
Since $\widehat{\theta}(m),~ \overline{\theta}(m),~ \dot{\theta}_i(m) \limitepsn \theta^*$ for $i=1,\cdots,|m|$, in this case of overfitting, the same arguments as in the proof of Proposition \ref{proposition} lead to
\begin{equation}\label{F_n_partial2_widehat_L_conv}
\widehat{F}_n(m) \limitepsn F(\theta^*,m)  \text{ and }  \dfrac{-2}{n} \dfrac{\partial^2 \widehat{L}_n\big( \overline{\theta}(m) \big)}{\partial \theta^2} \limitepsn F(\theta^*,m) ~  \text{ where } ~ F(\theta^*,m)= \Big(\E \Big[ \dfrac{ \partial^2 q_0(\theta^*)}{ \partial\theta_i \partial\theta_j}  \Big] \Big)_{i,j \in m} . 
\end{equation}     
For the overfitted model $m$, on can deduce from \cite{Bardet2009} that $F(\theta^*,m)$ is invertible, thus for $n$ large enough and with a sufficiently large probability, the matrix $\widehat{F}_n$ is invertible.   
 From (\ref{Taylor_widehat_L}), (\ref{theta_m_theta_star_partial_widehat_L}), (\ref{F_n_partial2_widehat_L_conv}) and Lemma \ref{lem1}, it holds that
 \begin{align}\label{L_m_star_L_m_log_log_n_partial}
 \nonumber &\frac{1}{\log \log n} \big(\widehat{L}_n\big(\widehat{\theta}(m)\big)-\widehat{L}_n\big(\widehat{\theta}(m^*)\big) \\
 \nonumber &=   \frac{2}{n^2 \log \log n}  \dfrac{\partial\widehat{L}_n(\theta^*)}{\partial \theta'}  \widehat{F}_n^{-1}(m) \dfrac{\partial^2 \widehat{L}_n\big( \overline{\theta}(m) \big)}{\partial \theta^2}  \widehat{F}_n^{-1}(m) \dfrac{\partial\widehat{L}_n(\theta^*)}{\partial \theta} \\
 \nonumber &=  - \Big( \frac{1}{ \sqrt{n \log \log n}}  \dfrac{\partial\widehat{L}_n(\theta^*)}{\partial \theta'} \Big) \widehat{F}_n^{-1}(m) \Big(\dfrac{-2}{n} \dfrac{\partial^2 \widehat{L}_n\big( \overline{\theta}(m) \big)}{\partial \theta^2} \Big) \widehat{F}_n^{-1}(m) \Big( \frac{1}{ \sqrt{n \log \log n}} \dfrac{\partial\widehat{L}_n(\theta^*)}{\partial \theta} \Big) \\
 &= \Big( \frac{1}{ \sqrt{n \log \log n}}  \dfrac{\partial L_n(\theta^*)}{\partial \theta'} + o(1) \Big) \mathcal{O}(1)  \Big( \frac{1}{ \sqrt{n \log \log n}}  \dfrac{\partial L_n(\theta^*)}{\partial \theta} + o(1) \Big) ~ ~ ~ a.s. ~.
 \end{align}
 Thus, according to (\ref{L_m_star_L_m_log_log_n_partial}) and (\ref{Proof_prop_LIL}), (\ref{L_m_star_L_m_log_log_n}) follows.\\
 Therefore, since $\kappa_n / \log \log n \limiten \infty$ and $|m| > |m^*|$, then (\ref{C_m_star_C_m_log_log_n}) and (\ref{L_m_star_L_m_log_log_n}) lead to
\begin{equation}\label{C_m_star_C_m_log_log_n_conv_infty}
 \lim_{n \rightarrow \infty} \frac{1}{ \log \log n} \big( \widehat{C}(m^*) - \widehat{C}(m) \big) =-\infty ~ ~ ~ a.s. ~ . 
\end{equation} 
This implies, 
\begin{equation}\label{C_m_star_C_m_n_large}
   \widehat{C}(m) - \widehat{C}(m^*) > 0  ~ ~ ~ a.s. ~ \text{ for large } n. 
\end{equation} 

\medskip

2. Let $m\in {\cal M}$ such as $m  \nsupseteq  m^*$. 
  We have,
  \begin{equation}\label{C_m_star_C_m_n}
 \frac{1}{  n} \big( \widehat{C}(m^*) - \widehat{C}(m) \big) = \frac{2}{ n} \big(\widehat{L}_n\big(\widehat{\theta}(m)\big)-\widehat{L}_n\big(\widehat{\theta}(m^*)\big) - \frac{\kappa_n}{ n}  ( |m|-|m^*|). 
\end{equation}
For all $\theta \in \Theta$, denote $L(\theta) = -\dfrac{1}{2} \E[q_0(\theta)]$.
According to the proof of Theorem 3.1 of \cite{Bardet2020}, we get
\begin{equation*} 
  \frac{1}{ n} \big(\widehat{L}_n\big(\widehat{\theta}(m)\big)-\widehat{L}_n\big(\widehat{\theta}(m^*)\big) = L(\theta^*(m)) - L(\theta^*) + o(1) ~ ~ a.s. ~ . 
 \end{equation*}
Note that, from \cite{Bardet2009}, the function $L: \Theta \rightarrow \R$ has a unique maximum at $\theta^*$.
Since $m  \nsupseteq  m^*$, it holds that $\theta^* \notin \Theta(m)$.
Hence, $L(\theta^*(m)) - L(\theta^*) < 0 ~ ~ a.s.$ .
Thus, according to (\ref{C_m_star_C_m_n}) and since $\kappa_n /n \limiten 0$, we get
\begin{equation}\label{C_m_star_C_m_n_conv_neg}
 \lim_{n \rightarrow \infty} \frac{1}{ n} \big( \widehat{C}(m^*) - \widehat{C}(m) \big) < 0 ~ ~ ~ a.s. ~ \text{ and } ~ ~  \widehat{C}(m)-\widehat{C}(m^*) >0  ~ ~ ~ a.s. ~ \text{ for large } n. 
\end{equation} 

\medskip
 
  Thus, the strong consistency of $ \widehat{m}=  \underset{m \in \mathcal{M}}{\text{argmin}} ~ \widehat{C}(m) = \underset{m \in \mathcal{M}}{\text{argmin}} \big( \widehat{C}(m) - \widehat{C}(m^*) \big)$ follows from (\ref{C_m_star_C_m_n_large}) and (\ref{C_m_star_C_m_n_conv_neg}). 
 \begin{flushright}
$\blacksquare$ 
\end{flushright}

\paragraph{Proof of Corollary \ref{corollary}}~ 

\medskip

According to the proof of Theorem \ref{theorem} (equations (\ref{C_m_star_C_m_n_large}) and (\ref{C_m_star_C_m_n_conv_neg})) it holds that $\widehat{m} = m^* ~ a.s.$ for large $n$. Thus, the corollary follows from Proposition \ref{proposition}. 
 \begin{flushright}
$\blacksquare$ 
\end{flushright}

 \bibliographystyle{acm}

  \end{document}